# Hanslick's Musical Formalism as a Mathematical Structure

*Examining meaning in Music Through the Lens of Category Theory*

Ilyas Khan[1]
734062@soas.ac.uk

We investigate Eduard Hanslick's formalist aesthetics through the prism of Category Theory. In his 1854 essay "On the Musically Beautiful" (*Vom Musikalisch-Schönen)*, Hanslick argued that music lacks external semantic reference, and defined musical content as "sonically (or tonally) moved forms." Counter intuitively, he also asserted that music possesses an internal, untranslatable "logic and sense." This article examines how Category Theory—the mathematics of universal structures and relationships—might provide a framework to better understand Hanslick's internal logic. Building on impressive advances in Applied Category Theory in the past few years, we propose an approach where modelling musical elements (pitches, rhythms, gestures) as objects and their transformations as morphisms, we illustrate how meaning in music can emerge relationally rather than referentially. We use the Yoneda Lemma to show how musical meaning is a function of structural compositionality and internal mapping. This approach can provide a rigorous, non-linguistic vocabulary to explain how absolute music can lack extramusical reference yet remain profoundly meaningful. We also hope that this work will help develop a far more nuanced and realistic appreciation of Hanslick than has been the case in modern musicological and ethnomusicological discourse.

## Part 1 Introduction and Background

The question of how music generates meaning, or if music has any meaning in the conventional sense that human beings think about when we refer to meaning in language, has been debated for almost as long as music has been discussed and analysed. Unlike day to day human language, which utilises a structured common lexicon of symbols and combinations of symbols that refer to external realities, and where semantic and syntactical rules provide a scaffolding that holds the overall structures together, absolute instrumental music (or 'pure' music), lacks clear semantic architectures. For example, a chord or note sequence does not refer to any external reality, such as "house" or "swan," nor does a descending minor third inherently mean "sadness".

Yet, at the same time, when listening to Brahm's Lullaby (in the case of music where the descending minor third is used), or a Bach fugue or a Shostakovich string quartet, we can have meaningful experiences and emotions that could include a sense of coherence, development, and even narrative purpose that exists in addition to and outside of any directed commentaries (such as titles for musical works, or concert notes by composers or publishers). This has led to conflicting views about whether, and how, music carries meaning, and how an art form that refers to nothing outside of itself (i.e. is a closed system), inspires us to project meaningfulness and emotion (if not exact meaning) onto it and from it.

*Vom Musikalisch-Schönen* (*On the Musically Beautiful*), written by Eduard Hanslick in the 1850's, was the first comprehensive treatment of many of these questions. For much of the Romantic era, pure music was romanticised and held up as an almost mystical conduit for the infinite, and a direct expression of unfiltered emotional states. In *On the Musically Beautiful* Hanslick argues very strongly against this prevailing "aesthetics of feeling" (*Gefühlsästhetik*). Hanslick argued[2] that contemporary musicology was fundamentally flawed because it confused the subjective impression a piece leaves on a listener with the objective properties of the work itself.

[1] SOAS University of London, and Topos Institute, Berkely CA

[2] On The Musically Beautiful by Eduard Hanslick (trans Payzant) Hackett,1986

Hanslick's article pitted him squarely against the prevalent views of society at that time, that was heavily influenced by the narrative of the Romantics, and which positioned music as a vessel for representing specific emotional states or external narratives[3]. Hanslick developed his basic thesis, that the beauty and essence of music are purely formal, in *On the Musically Beautiful*, but then went much further, and developed the argument that the content of music is nothing other than "sonically moved forms" (*tönend bewegte Formen*)[4]. For Hanslick, music does not *contain* meaning; it *is* the meaning, and music does not bring meaning to humans, but humans take meaning to music.

Hanslick's position that music is self-referential has been criticised by a wide spectrum of critics and philosophers who characterise him as the ultimate 'formalist[5]'. The key Hanslick opinion that caused so much anguish amongst his critics (that continues to modern times) is the view that music is not a mere acoustic amusement. In his essay he states:

> *"...in music there is sense and logic, but musical sense and logic. It is a language that we speak and understand, but we are unable to translate……The true heart of music, the formal beauty which gratifies in itself, would thereby be pierced through, (were the) 'chimera' of meaning be pursued[6]"*

This then leads to the obvious question "if music possesses an internal 'logic and sense' that is an untranslatable language, how can we study, describe, or validate this logic without relying on the very verbal metaphors that exist externally (i.e. external to the sound itself) that Hanslick's views have prohibited?

For over a century, traditional musicology has struggled to bridge this gap, with arguments for and against Hanslick ranging between very spare reductionist formal analysis on one extreme, and subjective, metaphor-heavy (and often emotion heavy) hermeneutics on the other. Even supporters of Hanslick have found this seeming dichotomy to be problematic[7].

It is not the purpose of this paper to defend (or attack) Hanslick. The purpose of this paper is to show how category theory (through its fundamental theorem, the Yoneda Lemma) brings clarity to the view that music itself (meaning the sounds) have no meaning other than with reference to themselves.

It is however, worth noting that Hanslick's formalism is frequently misunderstood as a cold, reductionist denial of musical value or beauty. As I note below, this could not be further from the truth, and is a mis-characterisation that, for example, Carl Dahlhaus notes, has predominated in formalist critiques for decades[8]. Hanslick did not view music as esoteric amusement or an assembly of empty acoustic patterns. On the contrary, he insisted that music is an intellectual and spiritual force of beauty. His thesis is, after all, titled “On the Musically Beautiful”.

---

[3] An excellent overview of this historical context is provided in Musicologica Austriaca in the entry by Landerer and Wilfing titled "Towards a Dynamic view of Hanslick's Aesthetics"

[4] Rethinking Hanslick ed Grimes,Donovan and Marx University of Rochester Press (2015) intro pp 1-9

[5] See for example Thomas Grey, "Hanslick, Eduard," in The New Grove Dictionary of Music and Musicians, ed. Stanley Sadie, 2nd ed. (London: Macmillan, 2001), 10:827–33. Grey summarises widespread negative views against Hanslick.

[6] On The Musically Beautiful by Eduard Hanslick (trans Payzant) Hackett,1986 Ch III

[7] See for example "What was Hanslick Denying ?" by Peter Kivy

[8] The Idea of Absolute Music, Carl Dahlhaus University of Chicago Press (1991)

Before we turn, therefore, to the core focus of this short article and look at how the mathematics of category theory can help reconcile the contradictions noted above, it is worth using Hanslick’s own words to frame the challenge that this paper hopes to at least partially address.

*"In music there is sense and logic, but musical sense and logic. It is a language that we speak and understand, but we are unable to translate[9]"*

This quote captures the paradox of Hanslick's aesthetics when taken entirely in isolation (meaning without reading the entirety of his essay and relying only on certain passages that might appear polemical). Hanslick is stating, above, that pure music is a medium that functions in many ways like human language, for example, possessing a highly developed, self-contained syntax that listeners intuitively comprehend, and yet remains completely untranslatable into spoken words or conceptual thoughts. The "sense" (*Sinn*) of a musical phrase is entirely intrinsic, according to Hanslick, and it cannot be exported outside the sonic boundaries of the piece being analysed or listened to, without losing all sense.

Traditional musicology has long struggled to resolve this tension. Linguistic structuralism and semiotics, which rely on the relationship between a signifier (the sound) and the signified (the external concept), fail when applied to absolute music because the "signified" in music is simply more musical structure. This leaves formalists with a deep methodological dilemma in how to rigorously analyze and validate internal, non-referential "logic" without defaulting to the verbal, emotional metaphors that Hanslick strictly prohibited.

This dilemma is more challenging when Hanslick also underscores his deeply and profoundly held views on the beauty of music. In the quote below, taken from the first edition of *On the Musically Beautiful*, but confusingly absent from many subsequent editions and translations, Hanslick sets out his stall as clearly as in any part of his monograph[10].

*“This spiritual content thus combines, in the soul of the listener, the beautiful in music with all the other great and beautiful ideas. He does not experience music merely as bare and absolute through its own beauty, but simultaneously as a sounding image of the great movements in the universe. Through deep and secret relationships to nature, the meaning of tones is heightened far beyond the tones themselves and allows always to feel the infinite even as we listen to the work of human talent. Just as the elements of music, - sound tone, rhythm, strength, weakness – are found in the entire universe, so man rediscover in music the entire universe”*

This quote does not alter the fact that Hanslick’s central thesis rests on a strict distinction between representation and form. When he argues that music is incapable of expressing definite emotions such as love, rage, or sorrow, he is doing so because any description of emotions requires clear cognitive concept, meaning at minimum, a framework of ideas, memories, social norms and judgments, which musical notes, intervals, and rhythms do not and cannot contain. As noted above, a descending minor third possesses no intrinsic semantic capability to denote "sadness" in the way the word *sadness* denotes a psychological state. If music cannot represent external psychological realities or worldly entities, what then is its content?

[9] On The Musically Beautiful by Eduard Hanslick (trans Payzant) Hackett,1986 p30
[10] Rethinking Hanslick ed Grimes, Donovan and Marx Rochester University Press (2015) p2

Hanslick's answer to this question forms the bedrock of his musical formalism, namely that the content and object of music are "sonically moved forms" (*tönend bewegte Formen*)[11]. Using the mathematics of category theory, it is this concept that we hope to clarify and at least partially address in this brief paper, It is precisely at this philosophical impasse that standard linguistic frameworks break down, and where a relational, category-theoretic ontology is a useful way of understanding meaning in music.

## Part 2a - Category Theory for Music

*"Since Greek antiquity it has been a tradition in European thinking to describe musical facts in a mathematical language" (Guerino Mazzola[12])*

This paper explores a possible resolution to Hanslick's paradox by introducing the structural mathematics of Category Theory to musical aesthetics. Developed by Samuel Eilenberg and Saunders Mac Lane in the mid-20th century[13], category theory is the mathematical study of universal structures and relationships. Rather than analysing mathematical entities by their internal composition and potentially falling foul of problems associated with reductionism, category theory defines mathematical entities (and as we shall see, other related entities) entirely by how they relate to, map into, and transform into one another. In recent years applied category theory has made very significant advances, not only in the sense that increasing numbers of mathematicians think of category theory as the meta language of mathematics, but also in applied mathematics, solving problems in real life, for example in areas such as software programming and scientific discovery[14].

As noted in the introduction, Hanslick puts forward the view that music possesses an internal, but untranslatable logic. The preliminary work contained in this paper examines this paradox through a formal framework that prioritises structure and relationship over external reference. category theory, the language of mathematics[15] provides precisely a vocabulary and framework that is required. As briefly noted above, in traditional mathematics, objects are typically defined by their internal contents (for instance, a set is defined by its elements). In category theory, however, an object is defined entirely by its structural relationships and pathways to other objects. This shifts the analytical focus from *what an entity is* to *how it behaves* within a larger ecosystem—a paradigm that perfectly mirrors the functional nature of musical elements, which derive their meaning entirely from their context within a work.

By modelling musical elements such as motifs, harmonic progressions, and rhythms as objects within a category, and their transformations as morphisms (or maps or transformations), category theory provides a precise, non-linguistic vocabulary required to define Hanslick's "internal logic" on a rigorous basis. Through advanced categorical concepts

---

[11] ibid p29

[12] The Topos of Music G Mazzola, Birkhauser (2002)

[13] The General Theory of Natural Equivalence" Mac Lane and Eilenberg (1945) and also Categories for the Working Mathematician Mac Lane and Eilenberg, Springer (1971)

[14] An Invitation to Applied Category Theory by B Fong and D Spivak Cambridge University Press (2019)

[15] See, for example, S Roman "The Language of Category Theory" Birkhauser (2017) and also E Riehl "Category Theory in Context" Dover (2016)

such as the Yoneda Lemma, this article examines ways in which we might mathematically demonstrate how meaning in absolute music emerges relationally and compositionally rather than referentially.

In this interdisciplinary framework this paper hopes to show that absolute music can be completely devoid of extramusical reference while remaining an intricately structured, profoundly logical, and inherently meaningful system.

***Section 2 is designed for readers with some very basic understanding of category theory. However they should be relatively easily accessible for the reader who not trained in maths. For that reader, or those not familiar with Category Theory and who wish to gain some basic knowledge that is more than intuition, we recommend Andre Cabannes' excellent book "Introduction to Category Theory[16]" that has been written for the average educated non-mathematician.***

In order to provide a rigorous underpinning for Hanslick's view that music has meaning even though that meaning is untranslatable, we propose a category that we name as the Category of Musical Structures (denoted CatMS), consisting of two foundational components, objects and morphisms (or arrows) that connect them. The objects (for example A, B, C… etc), represent distinct musical entities. The objects A, B, C can be a chord, a sequence of chords, a single pitch, a rhythmic group, a melodic motif, or an entire symphonic movement. Morphisms, or arrows, from one object to another, such as from A to B, represent the transformations, operations, or transitions that map or transform, one musical object into another.

CatMS complies with the axioms of category theory by fulfilling two criteria:

- **Compositionality:** If there is an arrow or morphism denoted F from Object A to Object B, (for example transposing a motif up a perfect fifth) and a subsequent transformation (or arrow) from Object B to Object C, (for example, inverting that transposed motif), there must exist a composite morphism or arrow from Object A to Object C, that represents the entire transformation
- **Associativity and Identity:** Associativity simply means that the order in which consecutive arrows or morphisms are grouped do not change the outcome. This means that if we have arrows from A to B, and then from B to C, we can denote them as either (A to B) and then to C, or A to (B and then to C). Identity means that every musical object must have an identity morphism that can be thought of as an action that does nothing, or an act of repetition that leaves the musical entity (i.e. the Object) unaltered.

Hanslick's concept of "sonically moved forms" is thus captured by thinking of the musical objects (A, B, C etc) as "forms," while the morphisms or arrows, are the "movement."

In order to illustrate this idea we consider standard western tonal operations and their progression, and use the example of a minor triad and a major triad as objects, whilst P (parallel), R (relative), and L (leading) tone-exchanges are the outcomes of morphisms or arrows that map these chords between each other. The diagram below should give the reader who is not familiar with category theory a good sense of how the basic structure holds together.

---

[16] "Introduction to Category Theory – journey in abstraction" Andre Cabannes

**Dia 1**

```
      L (Leading-Tone Exchange)
  C major <----------------------> E minor
    |                               ^

    | P (Parallel)                  | R (Relative)
    v                               |
  C minor <------------------------+
```
*(diagram generated with the assistance of the AI app Gemini)*

In addition to thinking about relatively simple morphisms (or shifts if you like that word instead) in the nature of the sounds within CatMS, category theory allows us to move beyond isolated chord changes or scales, through the use of functors. Think of functors as morphisms or maps that preserve the relationships between entirely different domains or categories. A relatively simple way to think of this is to imagine a functor ***F*** that goes from CatMS to CatMP (the category of musical performances), thus mapping abstract structural objects from a musical score (ie. CatMS) to physical actions within a category of orchestral gestures, actions, and instruments (ie. CatMP).

The key point to be borne in mind here is that since functors preserve compositionality, they guarantee that the internal logic established by the composer in the score translates into the physical movements of the performer or member of an orchestra, thereby showing that music is no longer seen as a sequence of static, isolated acoustic events. Music is a dynamic, interconnected network of transformations (maps and morphisms). Formalising the score as an atlas or map of maps, of objects and morphisms, lays the mathematical foundation required to demonstrate Hanslick's conviction that "meaning" in music emerges internally and structurally, completely independent of external linguistic reference.

We now show how this 'validation' can be achieved through the Yoneda Lemma.

### Part 2b - Internal Logic Through the Yoneda Lemma

"An object is completely determined by its relationships to other objects"

One of the most profound theorems in mathematics, and the fundamental theorem of category theory, is the Yoneda Lemma[17]. We believe that this theorem of category theory may help validate Eduard Hanslick's assertion that music possesses an internal "logic and sense" that, whilst intelligible or understandable in one sense to the listener(s), is rendered untranslatable or is incapable of being translated into words. This concept of something being understood but not being capable of being captured by language is not novel or rare, and was, for example, further developed in the 20th century by the early Wittgenstein and his ideas about the limits of our language[18].

---

[17] For a thorough description and exposition of the Yoneda Lemma see "Yoneda Lemma" by T Rud and J Ton, Ecole Polytechnique Federal de Lausanne (2013). A more basic, but in my view, very effective description is provided in "The Yoneda Lemma: Whats it all about?" by T Leinster, University of Edinburgh, Oct 2000

[18] "Die genzen meiner sprache bedeuten die grenzen meiner welt" or "The limits of my language are the limits of my world", Tractatus Logico-Philosophicus prop 5.6 (1921)

Before we continue forward into the Yoneda lemma and it's application to musicology, it is worth reminding ourselves that in conventional semiotic frameworks, a signifier derives its meaning via an extramusical semantic mapping in the form of a direct reference to an external object or psychological state. Hanslick argued that absolute music lacks this extramusical referent. This reliance of relationships between the different elements of sound existing within a musical work is the basis of our proposal that category theory, where the relationships between objects is what defines them, and the Yoneda lemma which formalises this position by rigorously showing how any object is defined only and absolutely by its relationships, can be so helpful.

The Yoneda Lemma therefore provides a uniquely interesting way to think about Hanslick's views[19].

The simplest statement of the Yoneda Lemma in mathematical terms (in this case provided by Tim Leinster from the University of Edinburgh and referenced in the footnotes below) is:

*[Cop ; Set], (Ha; X) is equivalent to X(A) where X is the functor acting on Object A*

For those not familiar with the notation, the statement simply explains that an object can be completely, uniquely, and exhaustively defined not by its internal essence or external reference, but by its relationships to all other objects within its category which, when applied to our example of music, can be explained to mean that any object A in the category CatMS, is entirely characterized by the functor that sends any other object in the category (think of these as Object X), to the set of morphisms Hom (X, A). This is the set where the elements are all possible maps and arrows from any object in CatMS to Object A. This concept is unpacked a little more carefully below.

Applying the Yoneda Lemma to CatMS, therefore implies that a musical entity, such as a single interval, a specific harmonic progression, or a complex melodic subject, possesses no intrinsic or static meaning in isolation. Instead, its identity and structural "sense" are mathematically equivalent to the sum total of all the ways other musical structures map into it, and how it maps into them[20].

**<u>Dia 2</u>**

```
              Hom(X, A)
 [Musical Object X] --------> [Musical Object A]
        |                          ^

        |                          |
        | α_X (Natural             |
        |     Transformation)      | Yoneda Embedding
        v                          |
  [Relational Network] ---------------+
```
*(diagram generated with the assistance of AI app Gemini)*

[19] Ch 9 from Topos of Music, G Mazzola provides a detailed perspective on how the Yoneda Lemma can be applied to the analysis of arts in general, and music in particular

[20] See appendix A for a more formal statement of the Yoneda Lemma

Diagram 2 shows how, through the Yoneda embedding, any attempt to analyse Object A, is identical to analyzing the entirety of relationships that A maintains with the rest of the objects in category CatMS. Since music contains a syntax (logical and structured 'scaffolding') and is a dense network of internal morphisms (transpositions, inversions, augmentations, and retrogrades), that are untranslatable into spoken language, they lend themselves to rigorous examination through category theory. When Hanslick posits that "meaning" is non-referential and is purely intentional and relational, he leads us to the logical conclusion that is captured by category theory. Just as any attempt to convert musical sounds into verbal definitions lacks 'sense' any attempt to translate a musical motif into a verbal word is to strip away its categorical network, destroying the very relationships that constitute its identity.

## Section 2c – An example of a simple Yoneda Embedding in CatMS the form of a Bach Fuge

By way of a reminder, our core category, CatMS is one who's objects A, B, C…. represent musical structures (e.g., pitch-class sets, motives, or contrapuntal matrices), and where, for any two objects A and B, the morphism (or arrow) from object A to object B consists of the structural transformations (for example retrograde-inversions, or rhythmic augmentations) mapping A to B.

We introduce the category CatSet (the category of sets), and functors from the category CatMSop (simply CatMS but with all the arrows between objects reversed) to CatSet. We further note the introduction of natural transformations. A natural transformation is simply a morphism or arrow between two functors.

The Yoneda Lemma shows that any natural transformation between a functor (call it X) applied to the musical object A in CatMSop and a separate functor that takes all the arrows or morphisms in CatMSop as they relate to object A, is the equivalent of the functor X acting on object A. Since that relationship (of a musical object) is intentional when applied to the set of all transformations, the notes or sounds remain defined only by the notes and sounds that they interact with, thus fulfilling the axiom laid down by Hanslick.

We provide a relatively simple example of how the Yoneda Lemma could work by looking at its application to *Johann Sebastian Bach's Fugue No. 2 in C minor (BWV 847*) from *The Well-Tempered Clavier*[21]. A fugue is an ideal musical entity for categorical analysis, as the whole of it's structural narrative is generated from a single sonic entity known as the subject.

In this example the five-measure fugue subject is represented as object A, (if its easier for the reader with little familiarity with musical terminology, you can also think of A as an arbitrary sequence of intervals). An Object B could be the alto voice that enters in measure 3 of the work, and where the transformation or morphism is a pitch transposition. Other objects are different parts of the work so that as the fugue progresses and Bach executes a sequence of transformations such as modal shifts and stretto fragmentations, the category is populated and defined (we remain in the category CatMS). Each of these transformations can be acted upon by functors that take these objects or sequences of objects, to the category CatSet, where the sum of all transformations to Object A (technically from Object A) can be placed since this sum of all transformations is a set.

[21] A free to view copy of the score is available online here

The "meaning" of the subject, for example at measure 15, or measure 20, is the historical and structural accumulation of its morphisms and the listener's aesthetic experience of the subject at the piece's climax is not a reaction to the acoustic properties of object A in isolation. Rather, without stretching the analogy too far, the listener cognitively tracks the natural transformations between functors acting from CatMSop (recalling that CatMSop is simply CatMS with the direction of all the arrows or morphisms reversed), to CatSet.

This reflects the reality that when a listener experiences the climax of the fugue, the aesthetic impact of the music is not derived from the notes themselves, but from the all too human processing of sound, and how the performance or recording in its entirety (or at any given point) is the sum of all its transformations upto that point.

**Dia 3**

```
Bach BWV 847 Fugue Subject Structure:
Measure 1Tonic Outline --> Measure 2Sequenced Scale --> Measure 3: Harmonic Cadence

      |                               |                                |
      v                               v                                v
   Morphism:                       Morphism:
Morphism:
Tonal Answer (G minor)          Soprano Inversion                 Bass Aug
```

*(diagram generated with the assistance of AI app Gemini)*

In the simplified diagram above, the five-measure C-minor subject of BWV 847 can be seen as a specific sequence of pitches and rhythms. Since we choose to treat the subject (as classically defined for a fugue), as  object A within CatMS there is no need for any further structure.

As the fugue unfolds, Bach constructs meaning relationally, perfectly illustrating the Yoneda principle in a process that can be reviewed for key stages as follows:

- **The Tonal Answer:** In measure 3, the alto voice introduces the subject transposed to the dominant key of G minor. This transformation is a morphism (Object A to Object B) that recontextualizes the original object.
- **The Counter-Subject:** As the answer enters, the original voice moves into a counter-subject. The structural integrity of the fugue is now defined by the simultaneous relationship between the subject and its counter-subject.
- **Contrapuntal Permutations:** Throughout the piece, Bach subjects the theme to inversion (flipping the intervals upside down), fragmentation, and modal shifts (moving into E-flat major) in such a manner that morphisms follow the rules of compositionality and associativity.

By defining the subject through a functor that takes Object A to the set of all morphisms to Object A in the category CatSet the Yoneda Lemma illustrates (and in a narrow sense, proves) that the musical "sense" championed by Hanslick (i.e. there are no outside references for meaning) is a mathematical reality. The fugue is a self-contained universe where objects are defined by their trajectories through morphisms, and provide a profound, completely non-linguistic form of cognitive meaning. The subject has no linguistic meaning in extensional terms, and any emotion or meaning that we attach is brought by us to the music, and not by the music to us.

As subject A is isomorphic to its relational functor (the functor that takes morphisms from other objects to A and places them in an object in CatSet that is the set of all morphisms), the Yoneda embedding proves that the musical "sense" Hanslick championed is a definable, mathematical reality. The "meaning" of the subject at any given moment in the fugue is no more (and no less !) than the historical accumulation of its morphisms.

**Section 3 - A Path to Resolving Formalist Dilemmas**

The formalisation noted above for the Yoneda Lemma, and its illustrative application to a Bach Fugue, provides a possible mathematical resolution to Hanslick's paradox of untranslatability. Using the framework offered by category theory, the "meaning" of a musical object A is not an intrinsic, isolated semantic entity. Rather, it is its image under the Yoneda embedding (via the functor that acts upon the object and maps it to the category CatSet). To understand a motif therefore, is to compute the actions of such functors for all contextually relevant structures (we refer to them above as X), in the fugue. The musical whole therefore operates as a self-contained universe where objects are defined entirely by their transformations, validating a profound, completely non-linguistic mode of syntax and meaning.

A perennial critique of musical formalism, and by extension, any mathematical treatment of musical architectures, is that it strips the medium of its emotional vitality and human phenomenological reality. Anti-formalists, from 19th-century expressive theorists to modern cognitive musicologists, argue that reducing music to cold syntax fails to account for the enthusiasm and emotions experienced by listeners. If music is merely a self-contained web of transformations, they claim, why does a specific rhythmic or melodic resolution evoke a profound sense of happiness, catharsis, sadness or tension?

By synthesizing Eduard Hanslick's aesthetics with category theory, however, we can demonstrate that this critique rests on a false dichotomy between structural logic and emotional expressiveness.

We suggest that the main advantage of the category theory model, specifically through the Yoneda embedding, is its rejection of a static, reductionist view of musical form. If at one extreme Hanslick is criticised by modern musicologists for being too formal, then by the same token, traditional set-theoretic or reductive mathematical methods of musical analysis often treat music as specimens that may as well be dead fragments of pitch and rhythm isolated on a page. Such static formalism struggles to capture the dynamic, temporal nature of listening. A category-theoretic formalism is inherently dynamic. Morphisms are not stationary points, they are vectors of change, translation, and development. When Hanslick defined music as tönend bewegte formen ("sonically moved forms"), his emphasis was fundamentally on the *movement*. His approach and philosophy was anything but static.

By treating the "meaning" of a musical object A as the presheaf of its relations, category theory operationalises this dynamic and kinetic energy. Human emotional responses to absolute music are not triggered by the isolated sonic object itself, but by the cognitive tracking of these relational shifts across time. Qualities and feelings such as "tension," "surprise," or "resolution" are the psychological results of navigating dense morphisms and pathways. A sudden harmonic deception or an unresolved dominant ninth, for example, represents a disruption in the expected natural transformations between representable functors. The listener's brain is not calculating algebra, but it is mapping structural trajectories, and when a trajectory is broken or elegantly subverted, an affective charge is generated. This is the reality of listening to music and allows us to anticipate and appreciate beauty of the kind that Hanslick so passionately advocates in the quote on p3 above, by understanding that this beauty and any related 'meaning' is extra-linguistic, and is not a contradiction, but a natural consequence of our interaction with the beauty of music

Although this is not within the scope of this article, we would point out that the relational framework afforded by category theory extends beyond the boundaries of the printed score. Using functorial mappings, it is possible to model the interface between the music-structural category CatMS and the category of human cognitive or gestural responses[22]. As functors preserve compositionality, a coherent internal relationship between two motifs in a Bach fugue are mapped directly onto a coherent internal relationship between cognitive or somatic states in the listener and gestures and movement by the instrumentalist.

Category theory bridges the gap between structure and sentiment. It reveals that musical emotion is not some vague, mystical vapor poured into the notes from some mysterious source. Expressiveness is an emergent property of structural compositionality. By proving that form is entirely composed of dynamic relationships, the Yoneda lemma allows us to defend Hanslick against his harshest critics. It helps to show that music can be wholly autonomous, entirely non-referential, and mathematically precise, while remaining a profoundly moving, deeply human but non-verbal language of pure relation.

In conclusion, we hope that this paper goes some way towards showing how the Yoneda Lemma provides a rigorous mathematical model for formalist analysis, showing that a system can be entirely self-contained without needing external anchors. It validates the position that musical grammar does not require external references to be fully specified; and even if it leaves aside the psychological or cultural aspects of how listeners subjectively experience sound, it allows for the fact that externalities are the source of meaning, meaningfulness or emotion. We hope this work will help towards proving Hanslick correct within the bounds of structural formalism and address the seeming contradiction noted above, and is a far more nuanced and realistic appreciation of Hanslick than has been pushed in modern musicological and ethnomusicological discourse.

***Acknowledgements***

I wish to acknowledge the support and guidance of my doctoral supervisor Professor Ilana Webster-Kogen who has been inspirational through her depth of knowledge and her intellectual brilliance. I also wish to thank Brendan Fong and David Spivak from the Topos Institute for the paths that they have opened-up for me in applying category theory to disciplines and interests far beyond conventional mathematics and computer science. I also wish to acknowledge Professor Bob Coecke who was instrumental in my engagement with category theory in the past dozen years.

[22] Gestural Pre-Sheaves: From Yoneda to Sheaves, Arias-Valero and Lluis-Puebla, MusMat (2021)

# Appendix

The Yoneda Lemma – a more traditional formal definition:

Let **C** be a locally small category. Let **Set** be the category of sets and let F: C→Set be a covariant functor. For any object A ∈ **C**, let $h_a$ = HOMc (A, -) denote the covariant hom-functor **C** to **Set**

The Yoneda lemma states that:

1. **Bijection** : There is a bijection (isomorphism in **Set**) between the set of natural transformations from $h_a$ to F and the set F(A): Nat($h_a$, F) ≡ F(A)
2. **The Mapping:** This bijection is explicitly given by evaluating a natural transformation @: hofA implies F at the identity morphism idA in Homc (A, A) : Phi(@) = @ A(idA)
3. **Naturality:** This isomorphism is **natural** in both A (as an object of C) and F (as an object of the functor category Set}^C)

**Selected Bibliography**